\documentclass[12pt]{article}
\usepackage{amsthm,amssymb,amsmath,graphicx}
\usepackage[width=5.35in,height=8.3in]{geometry}

\newtheorem{theorem}{Theorem}
\newtheorem{proposition}[theorem]{Proposition}
\newtheorem{corollary}[theorem]{Corollary}
\newtheorem{lemma}[theorem]{Lemma}
\theoremstyle{definition}
\newtheorem{example}[theorem]{Example}
\newtheorem{definition}[theorem]{Definition}
\newtheorem{remark}[theorem]{Remark}

\author{Luis Felipe Tabera\footnote{The author is supported by the project
MTM2005-08690-C02-02 and a FPU research grant from the Spanish Ministerio de
Educaci\'on y Ciencia.}\footnote{This article has been accepted for publication in the Revista Matem\'atica Iberoamericana http://www.uam.es/matem/ibero/}}
\title{Tropical Resultants for Curves and Stable Intersection}

\begin{document}
\maketitle
\begin{abstract}
We introduce the notion of resultant of two planar curves in the tropical
geometry framework. We prove that the tropicalization of the algebraic resultant
can be used to compute the stable intersection of two tropical plane curves. It
is shown that, for two generic preimages of the curves to an algebraic
framework, their intersection projects exactly onto the stable intersection of
the curves. It is also given sufficient conditions for such a generality in
terms of the residual coefficients of the algebraic coefficients of defining
equations of the curves.
\end{abstract}

\emph{keywords}: tropical geometry, resultants, plane curves


\section{Introduction}
In the context of tropical geometry, it is well known that two tropical curves
may share an infinite number of intersection points without sharing a common
component. This problem is avoided with the notion of stable intersection,
\cite{rgst}. Given two curves, there is a well defined set of intersection
points that varies continuously under perturbations of the curves. This stable
intersection has very nice properties. For example, it verifies a tropical
version of Bernstein-Koushnirenko Theorem (cf. \cite{rgst}). An alternative way
of defining a finite intersection set is the following: given two tropical
curves $f$ and $g$, take two algebraic curves $\widetilde{f}$ and
$\widetilde{g}$ projecting onto the tropical curves. Then, the intersection of
the two algebraic curves $\widetilde{f}\cap\widetilde{g}$ will project into the
intersection of the tropical curves, $T(\widetilde{f}\cap
\widetilde{g})\subseteq f\cap g$. In general, the set $T(\widetilde{f}\cap
\widetilde{g})$ depends on the election of the curves $\widetilde{f}$ and
$\widetilde{g}$. We are proving that, if the coefficients of $\widetilde{f}$,
$\widetilde{g}$ are generic, then the algebraic intersection will project
exactly onto the stable intersection and there is a correspondence among the
multiplicities of the intersection points. Moreover, given two curves, we
compute residually dense sufficient conditions defining these genericity
conditions. The method works in any characteristic and it is essential in the
generalization of the geometric construction method of \cite{Pappus-trop} to the
non
linear case. This also provides a particular case of Bernstein-Koushnirenko
theorem for fields of positive characteristic. In order to prove this
relationship, we introduce the notion of tropical resultant, as the
tropicalization of the algebraic resultant.

The paper is structured as follows. In Section \ref{sec_notions}, we give a
brief description of the algebraic context we will work in. Next, in Section
\ref{intersec_estable} we recall the notion of stable intersection for plane
curves and provide a brief discussion about its properties. Then, it is
introduced the notion of tropical resultant for univariate polynomials (Section
\ref{sec_univariate}) and plane curves (Section \ref{sec_bivariate}). In Section
\ref{sec_stable_int}, we relate the stable intersection of tropical curves with
the resultant of the curves and the generic preimage under the tropicalization
map. We will provide conditions for the lifts (preimages) to be compatible with
the stable intersection and the correspondence of the multiplicities. Finally,
in Section \ref{sec_remarks} we present some comments and remarks about the
results.

\section{Some basic notions in Tropical Geometry}\label{sec_notions}
The algebraic context where the theory is developed is the following:

Let $\mathbb{K}$ be an algebraically closed field provided with a non trivial
rank one valuation $v:\mathbb{K}^*\rightarrow \Gamma$. Without loss of
generality, we may suppose that $\mathbb{Q}\subseteq \Gamma \subseteq
\mathbb{R}$ and that $v$ is onto $\Gamma$. We denote by $k$ the residual field
of $\mathbb{K}$ under the valuation. We will distinguish two main cases along
the paper: the case whether $char(k)=0$ (hence $char(\mathbb{K})=0$) and the
case wether $char(k)=p>0$. In this case, either $char(\mathbb{K})=p$
(equicharacteristic $p$) or $char(\mathbb{K})=0$ ($p$-adic case). It is also
assumed that we have fixed a multiplicative subgroup $\Gamma' \subseteq
\mathbb{K}^*$ that is isomorphic to $\Gamma$ by the valuation map. The element
$t^\gamma$ represents the element of $\Gamma'$ whose valuation is $\gamma$. Any
element $x$ of $\mathbb{K}^*$ can be uniquely written as $x=x_0t^\gamma$, where
$v(x_0)=0$. We will denote the principal coefficient of an element $x$ of
$\mathbb{K}$ by $Pc(x)= \overline{x}_0 \in k^*$ and $Pc(0)=0$. We denote the
principal term of an element by $Pt(x)=\overline{x}_0t^\gamma$. The principal
term of an element is only a notation, it is not, in general, an element of
$\mathbb{K}$. If $y$ is an element of $\mathbb{K}^*$, $Pt(x)=Pt(y)$ if and only
if $v(x)=v(y)< v(x-y)$.

The tropicalization map is minus the valuation, $T(x)=-v(x)$. The tropical
semiring $\mathbb{T}$ is the group $\Gamma$ with the operations of tropical
addition $``a+b"=\max\{a,b\}$ and tropical product $``ab"=a+b$. With these
operations, $T(ab)=``T(a)T(b)"$ and, if $v(a)\neq v(b)$ or $v(a)=v(b)=v(a+b)$
then $T(a+b)=``T(a)+T(b)"$. Let $f=`` \sum_{i\in I} a_ix^i"=\max_{i \in
I}\{a_i+ix\}\in \mathbb{T}[x_1, \ldots, x_n]$ be a polynomial of support $I$,
where $x=x_1,\ldots,x_n$, $i=i_1,\ldots, i_n$, $ix=i_1x_1+ \ldots +i_nx_n$. The
set $\mathcal{T}(f)$ of zeroes of $f$ is the set of points in $\mathbb{T}^n$
such that the maximum of the piecewise affine function $\max_{i \in
I}\{a_i+ix\}$ is attained for at least two different indices. It is known
(Kapranov's Theorem, \cite{einsiedler-2004-}) that if $\widetilde{f}=\sum_{i \in
I}
\widetilde{a}_i x^i$ is any polynomial in $\mathbb{K}[x_1,\ldots, x_n]$ such
that $T(\widetilde{a}_i)=a_i$, then $T(\{\widetilde{f}(x)=0\}\cap
(\mathbb{K}^*)^n)$ is exactly the set of zeroes of $f$. Moreover, if $q\in
\mathbb{T}^n$ is a point, let $J\subseteq I$ be the set of indices where the
value $f(q)$ is attained and let $\alpha_i=Pc(\widetilde{a}_i)$. We define the
residual polynomial of $\widetilde{f}$ over $q$ as:
\begin{displaymath}
\widetilde{f}_q(x_1, \ldots, x_n) = \sum_{i \in J} \alpha_i x^i = Pc(
\widetilde{f} (x_1 t^{-q_1}, \ldots, x_n t^{-q_n})) \in k[x_1, \ldots, x_n]
\end{displaymath}

Then, it happens that:

\begin{theorem}\label{initialserie}
Let $\widetilde{f}\in \mathbb{K}[x_1, \ldots, x_n ]$ and $(\widetilde{b}_1,
\ldots, \widetilde{b}_n) \in (\mathbb{K}^*)^n$ be any point, then there is a
root $(\widetilde{c}_1, \ldots, \widetilde{c}_n)$ of $\widetilde{f}$ such that
$Pt(\widetilde{c}_i)=Pt(\widetilde{b}_i)$, $1\leq i\leq n$, if and only if
$b=(T(\widetilde{b}_1), \ldots, T(\widetilde{b}_n))$ is a zero of the tropical
polynomial $f$ and $(Pc(\widetilde{b}_1),\ldots, Pc(\widetilde{b}_n))$ is a root
of $\widetilde{f}_b$ in $(k^*)^n$.
\end{theorem}

For a constructive proof of this theorem we refer to \cite{Kapranov-EACA}
or \cite{Lifting-Constr}.

Let $C$ be a tropical plane curve defined as the zero set of a tropical
polynomial $f=``\sum_{(i,j)\in I}a_{(i,j)}x^iy^j"$. If we multiply $f$ by a
monomial, the curve it defines stays invariant. We define the support of $C$ as
the support $I$ of $f$ modulo a translation of an integer vector in
$\mathbb{Z}^2$. Analogously, given an algebraic curve $C$ in $(\mathbb{K}^*)^2$
defined by an algebraic polynomial $\widetilde{f}=\sum_{(i,j)\in
I}\widetilde{a}_{(i,j)}x^iy^j$, multiplying by a monomial does not change the
set of zeroes in the algebraic torus $(\mathbb{K}^*)^2$, we also define the
support of $C$ as the set $I$ modulo translations by an integer vector. If $C$
is a tropical curve, it may happen that there are polynomials with different
support defining $C$, even under the identification by translations we have
defined. Hence, when we define a tropical curve, we will always fix the support
of a defining polynomial.

Let $I$ be the support of a tropical polynomial $f$, the convex hull
$\Delta=\Delta(I)$ of $I$ in $\mathbb{R}^n$ is the Newton polytope of $f$. This
object is strongly connected with the set of zeroes of $f$. Every tropical
polynomial $f$ defines a regular subdivision of its Newton polytope $\Delta$.
The topological closure of $\mathcal{T}(f)$ in $\mathbb{R}^n$ has naturally a
structure of piecewise affine polyhedral complex. This complex is dual to the
subdivision induced to $\Delta$. To achieve this duality we have first to define
the subdivision of $\Delta$.

Let $\Delta'$ be the convex hull of the set $\{(i,t)| i\in I, t\leq
a_i\}\subseteq \mathbb{R}^{n+1}$. The upper convex hull of $\Delta'$, that is,
the set of boundary maximal cells whose outgoing normal vector has its last
coordinate positive, projects onto $\Delta$ by deleting the last coordinate.
This projection defines the regular subdivision of $\Delta$ associated to $f$
(See \cite{Mik05} for the details).

\begin{proposition}\label{Newton-Dual}
The subdivision of $\Delta$ associated to $f$ is dual to the set of zeroes of
$f$. There is a bijection between the cells of $Subdiv(\Delta)$ and the cells of
$\mathcal{T}(f)$ such that:
\begin{itemize}
\item Every $k$-dimensional cell $\Lambda$ of $\Delta$ corresponds to a cell
$V^\Lambda$ of $\mathcal{T}(f)$ of dimension $n-k$ such that the affine linear
space generated by $V^\Lambda$ is orthogonal to $\Lambda$. (In the case where
$k=0$, the corresponding dual cell is a connected component of $\mathbb{R}^n
\setminus \overline{\mathcal{T}(f)}$)
\item If $\Lambda_1\neq \Lambda_2$, then $V^{\Lambda_1}\cap V^{\Lambda_2} =
\emptyset$
\item If $\Lambda_1\subset \overline{\Lambda}_2$, then $V^{\Lambda_2}\subset
\overline{V^{\Lambda_1}}$
\item $\displaystyle{\mathcal{T}(f)=\bigcup_{0\neq\dim(\Lambda)}V^\Lambda}$
where the union is disjoint.
\item $V^\Lambda$ is not bounded if and only if $\Lambda\subseteq \partial
\Delta$.
\end{itemize}
\end{proposition}

From this, we deduce that, given a fixed support $I$, there are finitely many
combinatorial types of tropical curves with support $I$. These different types
are in bijection with the different regular subdivisions of $\Delta$.

Finally, let $C$ be a tropical planar curve of support $I$ and Newton polygon
$\Delta$, let $\Lambda$ be a one-dimensional cell of the subdivision of $\Delta$
dual to $C$, then, the weight of the dual cell $V^\Lambda$ is defined as
$\#(\overline{\Lambda}\cap \mathbb{Z}^2)-1$, the integer length of the segment
$\Lambda$.

\section{The Notion of Stable Intersection}\label{intersec_estable}

One of the first problems encountered in tropical geometry is that the
projective geometry intuition is no longer valid. If we define a tropical line
as the set of zeroes of an affine polynomial $``ax+by+c"$, then two different
lines always intersect in at least one point. The problem is that sometimes they
intersect in more than one point. The usual answer to deal with this problem is
using the notion of stable intersection.

Let $C_f$, $C_g$ be the set of zeroes of two tropical polynomials $f$ and $g$
respectively. Let $P$ be the intersection of the curves, $P=C_f\cap C_g$. It is
possible that $P$ is not the image of an algebraic variety $\widetilde{P}$ by
the map $T$. We want to associate, to each $q\in P$ an intersection
multiplicity. We will follow the notions of \cite{rgst} and we will compare them
with the subdivisions of the associated Newton polygons of the curves in terms
of mixed volumes. See \cite{Stu02} to precise the comparison between mixed
volumes and intersection of algebraic curves.

Let $C_{fg}=C_f\cup C_g$. It is easy to check that the union of the two tropical
curves is the set of zeroes of the product $``fg"$. The Newton polygon
$\Delta_{fg}$ of $C_f\cup C_g$ is the Minkowski sum of $\Delta_f$ and
$\Delta_g$. That is:
\[\Delta_{fg}=\{x+y\ |\ x\in\Delta_f, y\in \Delta_g\}\]
The subdivision of $\Delta_{fg}$ dual to $C_{fg}$ is a subdivision induced by
the subdivisions of $\Delta_f$, $\Delta_g$. More concretely, let $q$ be a point
in $C_{fg}$, let $\{i_1, \ldots, i_n\}$ be the monomials of $f$ where $f(q)$ is
attained and let $\{j_1, \ldots, j_m\}$ be the monomials of $g$ where $g(q)$ is
attained. Then $n\geq 2$ or $m\geq 2$. The monomials where $``fg"$ attains it
maximum are $\{i_rj_s|\, 1\leq r\leq n\,, 1\leq s\leq m\}$. The Newton polygon
of these monomials is the Minkowski sum of the Newton polygons of $\{i_1,
\ldots, i_n\}$ and $\{j_1, \ldots, j_m\}$, each one of these Newton polygons is
the cell dual to the cell containing $q$ in $\Delta_{fg}$, $\Delta_f$ and
$\Delta_{g}$ respectively. This process covers every cell of dimension $1$ and
$2$ of $\Delta_{fg}$. The zero dimensional cells correspond to points $q$
belonging neither to $C_f$ nor to $C_g$. Let $i, j$ be the monomials of $f$ and
$g$ where the value at $q$ is attained. Then the monomial of $``fg"$ where
$(``fg")(q)$ is attained is $ij$. To sum up, every cell of $\Delta_{fg}$ is
naturally the Minkowski sum of a cell $u$ of $f$ and a cell $v$ of $g$. The
possible combination of dimensions $(\dim(u), \dim(v), \dim(u+v))$ are:

\begin{itemize}
\item $(0,0,0)$, these cells do not correspond to points of $C_{fg}$.
\item $(1,0,1)$, these are edges of $C_{fg}$ that correspond to a maximal
segment contained in an edge of $C_f$ that does not intersect $C_g$.
\item $(2,0,2)$, correspond to the vertices of $C_{fg}$ that are vertices of
$C_f$ that do not belong to $C_g$.
\item $(1,1,2)$, this combination defines a vertex of $C_{fg}$ which is the
unique intersection point of an edge of $C_f$ with an edge of $C_g$.
\item $(1,1,1)$ are the edges of $C_{fg}$ that are the infinite intersection of
an edge of $C_f$ and an edge of $C_g$.
\item $(1,2,2)$ corresponds with the vertices of $C_{fg}$ that are a vertex of
$C_g$ belonging to an edge of $C_f$.
\item $(2,2,2)$ This is a vertex of $C_{fg}$ which is a common vertex of $C_f$
and $C_g$.
\end{itemize}
and the obvious symmetric cases $(0,1,1)$, $(0,2,2)$ and $(2,1,2)$.

If the relative position of $C_f$, $C_g$ is generic, then $C_{fg}$ cannot
contain any cell of type $(1,1,1)$, $(1,2,2)$ and $(2,2,2)$. That is, the
intersection points $q$ of $C_f$ and $C_g$ are always the unique intersection
point of an edge of $C_f$ and an edge of $C_g$. This is the transversal case.
The definition of intersection multiplicity, as presented in \cite{rgst} for
these cells $(1, 1, 2)$ is the following:

\begin{definition}\label{mult_rgst}
Let $q$ be an intersection point of two tropical curves $C_f$ and $C_g$. Suppose
that $q$ is the unique intersection line of an edge $r$ of $C_f$ and an edge $s$
of $C_g$. Let $\overrightarrow{r}$ be the primitive vector in $\mathbb{Z}^2$ of
the support line of $r$. Let $\overrightarrow{s}$ be the corresponding primitive
vector of $s$. Let $u$ be the dual edge of $r$ in $\Delta_f$ and let $v$ be the
dual edge of $s$ in $\Delta_g$, we call $m_u=\#(\overline{u} \cap \mathbb{Z}^2)
-1$ and $m_v=\#(\overline{v} \cap \mathbb{Z}^2) -1$ the weight of the edges $r$
and $s$ respectively. The intersection multiplicity is
\[mult(q)=\left|m_um_v\begin{vmatrix} \overrightarrow{r_x} &
\overrightarrow{r_y} \\
\overrightarrow{s_x} & \overrightarrow{s_y} \end{vmatrix} \right|\]
the absolute value of the determinant of the primitive vectors times the weight
of the edges.

If the curves are not in a generic relative position, consider the curve
$C_{f}^v$ obtained by translation of $C_f$ by a vector $v$. If the length of $v$
is sufficiently small, $|v|<\epsilon$ (that is, it is an infinitesimal
translation), then every cell of $\Delta_{f'g}$ of type $(0,0,0)$, $(1,0,1)$,
$(2,0,2)$ and $(1,1,2)$ stays invariant. Furthermore, if the translation is
generic (for all but finitely many directions of $v$), the cells of type
$(1,1,1)$ are subdivided into cells of type $(0,0,0)$ and $(0,1,1)$. That is, if
two edges intersect in infinitely many points, after the translation, every
intersection point will disappear. If $q$ is an intersection point of $C_f$ and
$C_g$ corresponding to a cell of type $(2,1,2)$ or $(2,2,2)$ and the direction
of $v$ is generic, this cell is subdivided, after the perturbation, into
cells of type $(0,0,0)$, $(1,0,1)$, $(1,1,2)$, $(2,0,2)$. That is, no
intersection point is a vertex of $f$ or $g$. However, some transversal
intersection points appear instead (of type (1,1,2)) in a neighborhood of $q$.
The intersection multiplicity of $q$ is, in this case, the sum of the
intersection multiplicities of the transversal intersection points.
\end{definition}

Now we recall the notion of stable intersection of curves (See \cite{rgst}).

\begin{definition}\label{interseccion_estable_def}
Let $C_f$, $C_g$ be two tropical curves. Let $C_f^v$, $C_g^w$ be two small
generic translations of $C_f, C_g$ such that their intersection is finite. The
stable intersection $C_f\cap_{st}C_g$ of $C_f$ and $C_g$ is the limit set of
intersection points of the translated curves $\lim_{v,w \rightarrow 0} (C_f^v
\cap C_g^w)$.
\end{definition}

From the previous comments it is clear that

\begin{proposition}
Let $C_f$, $C_g$ be two tropical curves, then the stable intersection of $C_f$
and $C_g$ is the set of intersection points with positive multiplicity.
\end{proposition}

This stable intersection has very nice properties. From the definition, it
follows that it is continuous under small perturbations on the curves. Moreover,
it verifies a Berstein-Koushnirenko Theorem for tropical curves.

\begin{theorem}\label{Bernstein}
Let $C_f$, $C_g$ be two tropical curves of Newton polygons $\Delta_f$,
$\Delta_g$. Then the number of stable intersection points, counted with
multiplicity is the mixed volumes of the Newton polygons of the curves
\[\sum_{q\in C_f \cap_{st} C_g} m(q)= \mathcal{M}(\Delta_f,\Delta_g)=
vol(\Delta_f +\Delta_g) -vol(\Delta_f) -vol(\Delta_g)\]
\end{theorem}
\begin{proof}
See \cite{rgst}
\end{proof}

In particular, we have the following alternative definition of intersection
multiplicity for plane curves:

\begin{corollary}\label{def_multiplicidad_tropical}
Let $f$, $g$ be two tropical polynomials of Newton polygons $\Delta_f$,
$\Delta_g$ respectively. Let $q\in \mathcal{T}(f) \cap \mathcal{T}(g)$ be an
intersection point. Let $\Lambda_f$, $\Lambda_g$ be the cells of
$Subdiv(\Delta_f)$, $Subdiv(\Delta_g)$ dual to the cells in the curve containing
$q$ respectively, then, the tropical intersection multiplicity of $q$ is:
\[mult(q)=\mathcal{M}(\Lambda_f, \Lambda_g)= vol(\Lambda_f + \Lambda_g) - vol(
\Lambda_f) - vol( \Lambda_g).\]
\end{corollary}
\begin{proof}
From the classification of intersection points, $q$ is an intersection point of
multiplicity zero if and only if it belongs to a cell of type $(1,1,1)$ in
$C_{fg}$. In this case $\mathcal{M}(\Lambda_f, \Lambda_g)= vol(\Lambda_f +
\Lambda_g) - vol( \Lambda_f) - vol( \Lambda_g)=0$, because an edge has no area.
If $q$ is a stable intersection point, let $f=``\sum_{i\in
\Delta_f}a_ix^{i_1}y^{i_2}", g=``\sum_{j\in \Delta_g}b_jx^{j_1}y^{j_2}"$, let
$f_q=``\sum_{i \in \Lambda_f} a_ix^i",\quad g_q=``\sum_{j\in \Lambda_g} b_jx^j"$
be truncated polynomials. It follows from the definition that the intersection
multiplicity of $q$ only depends in the behaviour of the mixed cell
$\Lambda_f+\Lambda_g$ in the dual subdivision of $\Delta_{fg}$. That is, the
intersection multiplicity of $q$ as intersection of $C_f$ and $C_g$ equals the
intersection multiplicity of $q$ as an intersection point of $\mathcal{T}(f_q)$
and $\mathcal{T}(g_q)$. But, by construction, the unique stable intersection
point of $\mathcal{T}(f_q)$ and $\mathcal{T}(g_q)$ is $q$ itself. Hence, by
Theorem~\ref{Bernstein}, the intersection multiplicity of $q$ is
\[\mathcal{M}(\Lambda_f, \Lambda_g)= vol(\Lambda_f + \Lambda_g) - vol(
\Lambda_f) - vol( \Lambda_g).\]
\end{proof}

\section{Univariate Resultants}\label{sec_univariate}
Let us start with the notion of tropical resultant of two univariate
polynomials. In algebraic geometry, the resultant of two univariate polynomials
is a polynomial that solves the decision problem of determining if both
polynomials have a common root.

\begin{definition}
Let $\widetilde{f}= \sum_{i=0}^{n} a_i x^i$, $\widetilde{g}= \sum_{j=0}^m b_j
x^j \in \mathbb{K}[x]$, where $\mathbb{K}$ is an algebraically closed field. For
simplicity, we assume that $a_0a_nb_0b_m\neq 0$. Let $p$ be the characteristic
of $\mathbb{K}$. Then, there is a unique polynomial in
$\mathbb{Z}/(p\mathbb{Z})[a_i,b_j]$, up to a constant factor, called the
resultant, such that it vanishes if and only if $\widetilde{f}$ and
$\widetilde{g}$ have a common root.
\end{definition}

In the definition, it is asked the polynomials to be of effective degree $n$ and
$m$, this is in order to avoid the specialization problems that usually appear
when using resultants. But the polynomials are also asked to have order zero.
This restriction is demanded for convenience with tropicalization. Recall that
the intersection of the varieties with the coordinate hyperplanes is always
neglected. Hence, the definition of resultant will take this into account.
Moreover, as the polynomials are described by its support, the resultant will
not be defined by the degree of the polynomials, but by their support. This
approach will be convenient in the next section, when there will be provided a
notion of resultant for bivariate polynomials.

\begin{definition}
Let $I$, $J$ be two finite subsets of $\mathbb{N}$ of cardinality at least 2
such that $0\in I\cap J$. That is, the support of two polynomials that do not
have zero as a root. Let $R(I,J,\mathbb{K})$ be the resultant of two polynomials
with indeterminate coefficients, $f=\sum_{i\in I} a_ix^i$, $g=\sum_{j\in J}
b_jx^j$ over the field $\mathbb{K}$. \[R(I, J, \mathbb{K}) \in \mathbb{Z}/(p
\mathbb{Z})[a, b],\] (where $p$ is the characteristic of the field
$\mathbb{K})$. Let $R_t(I, J, \mathbb{K})$ be the tropicalization of
$R(I,J,\mathbb{K})$. This is a polynomial in $\mathbb{T}[a,b]$, which is called
the tropical resultant of supports $I$ and $J$ over $\mathbb{K}$.
\end{definition}

So, our approach is to define the tropical resultant polynomial as the
projection of the algebraic polynomial. In this point, one may obtain, for the
same support sets $I$ and $J$, different tropical resultants, one for each
possible characteristic of $\mathbb{K}$. This is not good, in the sense that
tropical geometry should not be determined by the characteristic of the field we
have used to define the projection. Hence, one has to take care of what is the
common information of these polynomials. The answer is complete: the tropical
variety they define is always the same. This variety is the image of any
resultant variety over a field $\mathbb{K}$, so it will code the pairs of
polynomials with fixed support that have a common root.

\begin{lemma}
The tropical variety $\mathcal{T}(R_t(I,J,\mathbb{K}))$ does not depend on the
field $\mathbb{K}$, but only on the sets $I$ and $J$.
\end{lemma}
\begin{proof}
Let $\mathcal{N}$ be the Newton polytope of the resultant defined over a field
$\mathbb{L}$ of characteristic zero, $\mathcal{N}\subseteq \mathbb{R}^{n+m+2}$.
It is known that the monomials of $R(I,J,\mathbb{L})$ corresponding to vertices
of $\mathcal{N}$ (extreme monomials) have always as coefficient $\pm 1$ (See,
for example, \cite{GKZ-polytope-resultant} or
\cite{Sturmfels-polytope_resultant}). Hence, the
extreme monomials in $R(I,J,\mathbb{K})$ are independent of the characteristic
of the field $\mathbb{K}$ and so is $\mathcal{N}$. If $x=(x_1^{j_1},\ldots,
x_N^{j_N})$ is a monomial of $R(I,J,\mathbb{K})$ that does not correspond to a
vertex of $\mathcal{N}$, then $x=\sum \lambda_i v_i$, $0\leq \lambda_i \leq 1$,
where $v_i=(x_{1}^{j_{i,1}}, \ldots, x_{N}^{j_{i,N}})$ are vertices of
$\mathcal{N}$. $T(\textrm{coeff} (v_i))= T(\pm 1) = 0$ and, as coeff($x$) is an
integer (or an integer mod $p$), it is contained in the valuation ring, that is,
$0 \geq T( \textrm{coeff}(x)) \in \mathbb{T} \cup \{-\infty\}$.
$T(\textrm{coeff}(x))$ is finite and not zero if and only if we are dealing with
a $p$-adic valuation and $p$ divides $\textrm{coeff}(x)$. It is $-\infty$ if and
only if the characteristic of $\mathbb{K}$ divides the coefficient. Hence, for
any evaluation $w$ of the indeterminates, we have that
$T(\textrm{coeff}(x))+w_1j_1+\ldots +w_N j_N\leq w_1j_1+\ldots +w_N j_N= \sum
\lambda_i v_i(w)\leq \max_i\{v_i(w)\}$, where $v_i(w)=w_1j_{i,1} + \ldots + w_N
j_{i,N}$. It follows that the maximum of the piecewise affine function
$R_t(I,J,\mathbb{K})$ is never attained in the monomial $x$ alone and that $x$
does not induce any subdivision in the cell it is contained. Thus, this monomial
does not add anything to the tropical variety defined by $R_t(I,J,\mathbb{K})$.
The tropical hypersurface $\mathcal{T}(R_t(I,J,\mathbb{K}))$ is, as a polyhedral
complex, dual to the subdivision of $\mathcal{N}$ induced by
$R_t(I,J,\mathbb{K})$ (cf. \cite{Mik05}). In this case, the subdivision of
$\mathcal{N}$ induced by the tropical polynomial is $\mathcal{N}$ itself. So
$\mathcal{T}(R_t (I, J, \mathbb{K}))$ is always the polyhedral complex dual to
$\mathcal{N}$ centered at the origin. This complex is independent of
$\mathbb{K}$.
\end{proof}

Hence, fixed two supports $I$, $J$, there may be different tropical polynomials
that can be called the resultant of polynomials of support $I$ and $J$. However,
the variety all of them define is always the same, so there is a good notion of
resultant variety. Now we prove that the resultant variety
$\mathcal{T}(R_t(I,J,\mathbb{K}))$ has the same geometric meaning than the
algebraic resultant variety.

\begin{lemma}
Let $I, J$ be two support subsets as before. Let $f=``\sum_{i\in I} a_ix^i"$,
$g=``\sum_{j\in J}^m b_jx^j"$ be two univariate tropical polynomials of support
$I$ and $J$. Then, $f$ and $g$ have a common tropical root if and only if the
point $(a_i,b_j)$ belongs to the variety defined by $R_t(I,J,\mathbb{K})$.
\end{lemma}
\begin{proof}
Suppose that $(a_i,b_j)$ belongs to $R_t(I,J,\mathbb{K})$. By
Theorem~\ref{initialserie}, we can compute an element $(\widetilde{a}_i,
\widetilde{b}_j)$ in the variety defined by $R(I,J,\mathbb{K})$. In this case,
$\widetilde{f}= \sum_{i \in I} \widetilde{a}_i x^i$ and $\widetilde{g}= \sum_{j
\in J} \widetilde{b}_j x^j$ are lifts of $f$ and $g$. That is,
$T(\widetilde{f})=f$, $T(\widetilde{g})=g$. Moreover, their coefficients belong
to the algebraic resultant, so the algebraic polynomials have a common root
$\widetilde{q}$ that is non zero by construction ($0\in I\cap J$). Projecting to
the tropical space, $f$ and $g$ have a common root $T(\widetilde{q})$.
Conversely, if $f$ and $g$ have a common root $q$, we may take any lift
$\widetilde{g}=\sum_{j\in J}\widetilde{b}_jx^j$ of $g$. Then, by
Theorem~\ref{initialserie}, we may lift $q$ to a root $\widetilde{q}$ of
$\widetilde{g}$. Finally, note that the coefficients of $f$ belong to the
hyperplane defined by the equation $\sum_{i\in I} z_i q^i$, so it can be lifted
to an algebraic solution $\widetilde{a}$ of the affine equation $\sum_{i\in I}
z_i \widetilde{q}^i$, the polynomial $\widetilde{f}=\sum_{i\in I}
\widetilde{a}_i x^i$ projects onto $f$ and has $\widetilde{q}$ as a root. By
construction, $\widetilde{f}$, $\widetilde{g}$ share a common root
$\widetilde{q}$, hence, their coefficients $(\widetilde{a}_i, \widetilde{b}_j)$
belong to the algebraic resultant variety. Projecting again, the coefficient
vector ($a_i$, $b_j$) of $f$ and $g$ belong to the tropical resultant.
\end{proof}

This Lemma about the geometric meaning of the resultant also shows that the
variety defined by $R_t(I,J,\mathbb{K})$ does not depend on the field
$\mathbb{K}$. At least as a set of points, because the tropical characterization
of two tropical polynomials having a common root does not depend on the field
$\mathbb{K}$.

\begin{example}
Consider the easiest nonlinear case, $I=J=\{0,1,2\}$, the resultant of two
quadratic polynomials. If $f=a+bx+cx^2$, $g=p+qx+rx^2$, the algebraic resultant
in characteristic zero is $R_0=r^2a^2 -2racp +c^2p^2 -qrba -qbcp +cq^2a +prb^2$
and, over a characteristic 2 field it is $R_2 = r^2a^2 +c^2p^2 +qrba +qbcp
+cq^2a +prb^2$. If $char(k)\neq 2$, the tropical polynomial is $P_1=``0r^2a^2
+0racp +0c^2p^2 +0qrba +0qbcp +0cq^2a +0prb^2"$. If $char(\mathbb{K})=0$ and
$char(k)=2$, the tropical polynomial is $P_2=``0r^2a^2 +(-1)racp +0c^2p^2 +0qrba
+0qbcp +0cq^2a +0prb^2"$. Finally, if $char(k)=char(\mathbb{K})=2$ then the
tropical polynomial is $P_3=``0r^2a^2 +0c^2p^2 +0qrba +0qbcp +0cq^2a +0prb^2"$.
The unique difference among these polynomials is the term $racp$. This monomial
lies in the convex hull of the monomials $r^2a^2$ and $c^2p^2$ and it does not
define a subdivision because its tropical coefficient is always $\leq 0$. The
piecewise affine functions $\max\{2r+2a, r+a+c+p, 2c+2p\}$, $\max\{2r+2a,-1
+r+a+c+p, 2c+2p\}$ and $\max\{2r+2a, 2c+ap\}$ are the same. So the three
polynomials define the same tropical variety.
\end{example}

\section{Resultant of Two Curves}\label{sec_bivariate}
In this section, the notion of univariate resultant is extended to the case
where the polynomials are bivariate.

\begin{definition}
Let $\widetilde{f}$ and $\widetilde{g}$ be two bivariate polynomials. In order
to compute the algebraic resultant with respect to $x$, we can rewrite them as
polynomials in $x$. \[\widetilde{f}=\sum_{i\in I}\widetilde{f}_i(y)x^i, \quad
\widetilde{g}=\sum_{j\in J}\widetilde{g}_j(y)x^j,\] where
\[\widetilde{f}_i=\sum_{k=o_i}^{n_i}A_{ik}t^{-\nu_{ik}}y^k,\
\widetilde{g}_j=\sum_{q=r_j}^{m_j}B_{jq}t^{-\eta_{jq}}y^q\] and $A_{ik}$,
$B_{jq}$ are elements of valuation zero. Let
$P(a_i,b_j,\mathbb{K})=R(I,J,\mathbb{K})\in \mathbb{Z}/(p\mathbb{Z})[a_i,b_j]$
be the algebraic univariate resultant of supports $I$, $J$. The algebraic
resultant of $\widetilde{f}$ and $\widetilde{g}$ is the polynomial
$P(\widetilde{f}_i, \widetilde{g}_j, \mathbb{K})\in \mathbb{K}[y]$. Analogously,
let $f=T(\widetilde{f})$, $g=T(\widetilde{g})$, $f=``\sum_{i\in I}f_i(y)x^i"$,
$g=``\sum_{j\in J}g_j(y)x^j"$, where \[f_i= ``\sum_{k=o_i}^{n_i} \nu_{ik} y^k",\
g_j= ``\sum_{q=r_j}^{m_j} \eta_{jq} y^q".\] Let $P_t(a_i, b_j, \mathbb{K})
=R_t(I, J, \mathbb{K}) \in \mathbb{T} [a_i, b_j]$ be the tropical resultant of
supports $I$ and $J$. Then, the polynomial $P_t(f_i, g_j, \mathbb{K})\in
\mathbb{T}[y]$ is the tropical resultant of $f$ and $g$.
\end{definition}

Again, we have different tropical resultant polynomials, one for each possible
characteristic of the fields $\mathbb{K}$ and $k$. We want to check that this
notion of tropical resultant also has a geometric meaning. In the algebraic
setting, the roots of the resultant $P(\widetilde{f}_i, \widetilde{g}_j,
\mathbb{K})$ are the possible $y$-th values of the intersection points of the
curves defined by $\widetilde{f}$ and $\widetilde{g}$. This is not the case of
the tropical resultant, because $P_t(f_i,g_j, \mathbb{K})$ only has finitely
many tropical roots, while the intersection $\mathcal{T}(f)\cap \mathcal{T}(g)$
may have infinitely many points and there may be infinitely many possible values
of the $y$-th coordinates. Again, this indetermination is avoided with the
notion of stable intersection. We will prove that the roots of $P_t(f_i, g_j,
\mathbb{K})$ are the possible $y$-th values of the stable intersection
$\mathcal{T}(f)\cap_{st}\mathcal{T}(g)$. This will be made in several steps, the
first one is to check that $T(V(P(\widetilde{f}_i,\widetilde{g}_j,
\mathbb{K})))=\mathcal{T}(P_t(f_i, g_j, \mathbb{K}))$, provided that $A_{ik},
B_{jq}$ are residually generic. Sometimes, for technical reasons, it is better
to work with an affine representation of the polynomials. The set of polynomials
$\widetilde{f}=\sum_{i\in I} a_i x^{i^1}y^{i^2}$ of fixed support $I$ is an open
subspace (a torus) of a projective space. The projective coordinates of
$\widetilde{f}$ are its coefficients $[a_i: i\in I]$. We may fix an index
$i_0\in I$. Then, the affine representation of $\widetilde{f}$ with respect to
this index is obtained by setting $a_{i_0}=1$ and $a_{i_j}=a_{i_j}/a_{i_0}$. We
prove that this dehomogenization process is also compatible with
tropicalization. That is, if we divide each algebraic coefficient $A_{ik} t^{-
\nu_{ik}}$ and $B_{jq} t^{-\eta_{jq}}$ by $A_{i_0 k_0} t^{-\nu_{i_0 k_0}}$ and
$B_{j_0 q_0}t^{-\eta_{j_0 q_0}}$ respectively and substitute each coefficient
$\nu_{ik}$, $\eta_{jq}$ of the tropical polynomials $f$ and $g$ by $\nu_{ik} -
\nu_{i_0 k_0}=``\nu_{ik}/ \nu_{i_0 k_0}"$ and $\eta_{iq}-\eta_{i_0 q_0}$
respectively, still we have that $T(V(P(\widetilde{f}_i, \widetilde{g}_j,
\mathbb{K})))=\mathcal{T}(P_t(f_i, g_j, \mathbb{K}))$.

\begin{lemma}\label{tropicalizar_resultante_no_homogenea_char0p}
Let $\widetilde{f}= \sum_{i\in I} \widetilde{f}_ix^i$, $\widetilde{g}=
\sum_{j\in J} \widetilde{g}_j x^j\in \mathbb{K}[x, y]$, where the coefficients
are \[\widetilde{f}_i= \sum_{k=o_i}^{n_i} A_{ik}t^{-\nu_{ik}} y^k,\
\widetilde{g}_j =\sum_{q=r_j}^{m_j} B_{jq} t^{-\eta_{jq}}y^q\] and let $f=
\sum_{i\in I} f_i(y)x^i$, $g=\sum_{j\in J} g_j(y) x^j$,
\[f_i=``\sum_{k=o_i}^{n_i} \nu_{ik} y^k", g_j= ``\sum_{q=r_j}^{m_j} \eta_{jq}
y^q"\] be the corresponding tropical polynomials. Suppose that $A_{ik}$,
$B_{jq}$ are residually generic. Then $T(V(P(\widetilde{f}_i, \widetilde{g}_j,
\mathbb{K})))=\mathcal{T}(P_t(f_i, g_j, \mathbb{K}))$.
\end{lemma}
\begin{proof}
First, we suppose that $char(k)=0$. In general, the composition of polynomials
does not commute with tropicalization, because, in the algebraic case, there can
be a cancellation of terms when performing the substitution that does not occur
in the tropical case. Recall that, by the nature of tropical operations, a
cancellation of terms in the tropical development of the polynomial never
happens. So, we have to check that there is never a cancellation of terms in the
algebraic setting. First, it is proved that there is no cancellation of
monomials when substituting the variables by polynomials without dehomogenizing.
$P(a_i,b_j,\mathbb{K})$ is homogeneous in the set of variables $a_{i}$ and in
the set of variables $b_j$. As the substitution is linear in the variables
$A_{ik}$ and $B_{jq}$, $P(\widetilde{f}_i, \widetilde{g}_j, \mathbb{K})$ is
homogeneous in $A_{ij}$ and $B_{jq}$. If we have two different terms $T_1$,
$T_2$ of $P(a_i, b_j, \mathbb{K})$, then there is a variable with different
exponent in both terms. Assume for simplicity that this variable is $a_1$ with
degrees $d_1$ and $d_2$ respectively. After the substitution, the monomials
obtained by expansion of $T_1$ are homogeneous of degree $d_1$ in the set of
variables $A_{1k}$ and the monomials coming from $T_2$ are homogeneous of degree
$d_2$ in the variables $A_{1k}$. Thus, it is not possible to have a cancellation
of terms and we can conclude that the homogeneous polynomial projects onto the
tropical homogeneous polynomial.

In the case we dehomogenize $\widetilde{f}$ and $\widetilde{g}$ with respect to
the indices $(i_0 k_0)$, $(j_0 q_0)$ respectively. By the homogeneous case, we
can suppose that all the variables $a_i\neq a_{i_0}$ and $b_j\neq b_{j_0}$ in
$P(a_i, b_j, \mathbb{K})$ have already been substituted by the polynomials
$\widetilde{f}_i$ and $\widetilde{g}_j$ respectively. The only possibility to
have a cancellation of terms is if there are two monomials of the form
$Xa_{i_0}^{d_1}b_{j_0}^{d_2}$, $Xa_{i_0}^{d_3}b_{j_0}^{d_4}$ with
$d_1+d_2=d_3+d_4$ and $X$ is a monomial in the variables $A_{ik}$, $B_{jq}$.
But, as the polynomial is multihomogeneous in $A$ and $B$, it must happen that
$d_1=d_3$ and $d_2=d_4$. That is, the original monomials were the same. So, a
cancellation of terms is not possible and the dehomogenized polynomial projects
into the dehomogenized tropical polynomial. In particular,
$T(V(P(\widetilde{f}_i, \widetilde{g}_j, \mathbb{K})))=\mathcal{T}(P_t(f_i, g_j,
\mathbb{K}))$.

Now suppose that $char(k)=p>0$. In this case, it is not necessarily true that
the tropicalization of the algebraic resultant is the tropical resultant. But we
are going to check that the monomials where these two tropical polynomials
differ do not add anything to the tropical variety $\mathcal{T}(f_i, g_i,
\mathbb{K})$. So, we are going to compare the monomials in $P(\widetilde{f}_i,
\widetilde{g}_j, \mathbb{K})$ and $P_t(f_i, g_j, \mathbb{K})$. The support of
both polynomials is contained in the support of $P_t(f_i, g_j, \mathbb{L})$,
where $\mathbb{L}$ is an equicharacteristic zero field. The first potential
difference in the monomials are those obtained by expansion of a monomial $m$ of
the univariate resultant $P(a_i, b_j, \mathbb{K}) = R(I, J, \mathbb{K})$ whose
coefficient has valuation in $[-\infty,0)$. That is, $p$ divides
$\textrm{coeff}(m)$. It happens that $m$ is never a extreme monomial. That is,
$m=\sum_{l} \lambda_l v_l$, $0\leq \lambda_l \leq 1$ and $v_i$ are extreme
monomials. So, for every $r$, $\textrm{coeff}(m) +m(f_i(r), g_j(r)) \leq
m(f_i(r), g_j(r))= \sum_l \lambda_l v_l(f_i(r), g_j(r)) \leq \max \{v_l(f_i(r),
g_j(r)) \}$. Hence, the monomials of $m(f_i(y), g_j(y))$ never add anything to
the tropical variety defined by $P(\widetilde{f}_i, \widetilde{g}_j,
\mathbb{K})$, because they are never greater than the monomials that appear by
the extreme monomials. The other source of potential differences in the
monomials is the decreasing of the tropicalization of some terms of the power
$(\sum_{k=o_i}^{n_i} A_{ik} t^{-\nu_{ik}} y^k)^N$ due to some combinatorial
coefficient $\binom{N}{m}$ divisible by $p$. But, in the tropical context, it
happens that
\[(``\sum_{k=o_i}^{n_i}\nu_{ik}y^k")^N=``\sum_{k=o_i}^{n_i}\nu_{ik}^Ny^{kN}"\]
as piecewise affine functions. The rest of the terms in the expansion do not
contribute anything to the tropical variety. The only terms that may play a role
are $\nu_{ik}^N$, $\eta_{jq}^M$. So, even if the tropicalization of the
polynomials $P( I, J,\mathbb{K})$ depends on the algebraic field $\mathbb{K}$,
the tropical variety they define is always the same and it is the tropical
variety defined by $P_t(f_i, g_j, \mathbb{K})$, including the weight of the
cells.
\end{proof}

So, the previous Lemma provides a notion of tropical resultant for bivariate
polynomials with respect to one variable. They also prove that this polynomials
define the same variety as the projection of the algebraic resultant in the
generic case. Our next goal is to provide a geometric meaning to the roots of
the tropical resultant in terms of the stable intersection of the curves.

\section{Computation of the Stable Intersection}\label{sec_stable_int}

Let $f$ be a tropical polynomial of support $I$ defining a curve, let $\Delta_f$
be the convex hull of $I$. By Proposition~\ref{Newton-Dual}, the coefficients of
$f$ induce a regular subdivision in $\Delta_f$ dual to $f$. This subdivision is
essential in the definition of tropical multiplicity and stable intersection for
the case of curves. Next, it is proved that, for sufficiently generic lifts
$\widetilde{f}$ and $\widetilde{g}$, their intersection points correspond with
stable intersection points of $f$ and $g$.

\begin{lemma}\label{generic-lift_stable-int}
Let $f$ and $g$ be two tropical polynomials in two variables. Let $L$ be its
stable intersection. Then, for any two lifts $\widetilde{f}$, $\widetilde{g}$
such that their coefficients are residually generic, the intersection of the
algebraic curves projects into the stable intersection. \[T(\widetilde{f} \cap
\widetilde{g}) \subseteq \mathcal{T}(f) \cap_{{st}} \mathcal{T}(g)\]
\end{lemma}
\begin{proof}
If every intersection point of $f$ and $g$ is stable, then there is nothing to
prove. Let $q$ be a non stable intersection point. This means that $q$ belongs
to the relative interior of two parallel edges of $\mathcal{T}(f)$ and
$\mathcal{T}(g)$. The residual polynomials $\widetilde{f}_q$ and
$\widetilde{g}_q$ can be written (after multiplication by a suitable monomial)
as $\widetilde{f}_q= \sum_{i=0}^n \alpha_i (x^ry^s)^i$, $\widetilde{g}_q=
\sum_{j=0}^m \beta_i (x^ry^s)^j$. If $\widetilde{f}$, $\widetilde{g}$ have a
common point projecting into $q$ then there is an algebraic relation among their
residual coefficients. Namely, the resultant of the polynomials $\sum_{i=0}^n
\alpha_i z^i$, $\sum_{j=0}^m \beta_i z^j$ with respect to $z$ must vanish. If
the residual coefficients of $\widetilde{f}$, $\widetilde{g}$ do not belong to
the resultant defined by each non stable intersection cell, the intersection in
the torus of $\widetilde{f}$, $\widetilde{g}$ projects into the stable
intersection of $f$ and $g$.
\end{proof}

So, there is a natural relation between the stable intersection of two tropical
curves $f$ and $g$ and the intersection of two generic lifts $\widetilde{f}$ and
$\widetilde{g}$ of the curves. On the other hand, the intersection of two
generic lifts can be determined by the algebraic resultant of the defining
polynomials. Applying tropicalization, this relationship links the notion of
stable intersection with the resultants. To achieve a true bijection between the
roots of the resultant and the intersection points of the curves, it is used the
relationship between the tropical and algebraic resultants. So, one needs to
concrete the generality conditions for the values values $A_{ik}$, $B_{jq}$ that
makes Lemma~\ref{tropicalizar_resultante_no_homogenea_char0p} and
Proposition~\ref{tropicalizar_resultante_no_homogenea_char0p} hold. Next
Proposition shows how to compute the residually conditions for the compatibility
of the resultant.

\begin{proposition}\label{resultantes}
Let $\widetilde{f}$, $\widetilde{g}\in\mathbb{K}[x,y]$. Then, there is a finite
set of nonzero polynomials in the principal coefficients of the coefficients of
$\widetilde{f}$, $\widetilde{g}$, that depends only on the tropicalization $f$
and $g$ such that, if no one of them vanishes, then
\[T(\textrm{Res}_x(\widetilde{f},\widetilde{g}))=\mathcal{T}(R(I,J,\mathbb{K})(f
,g)).\] Where $R(I,J,\mathbb{K})(f,g)$ is the evaluation of the tropical
resultant of supports $I$ and $J$ in the coefficients of $f$ and $g$ as
polynomials over $x$.
\end{proposition}

\begin{proof}
Write $\widetilde{f}= ``\sum_{i,k} \widetilde{a}_{ik} x^i y^k"$,
$\widetilde{g}=`` \sum_{j,q} \widetilde{b}_{jq} x^j y^q"$, and take
$Pc(\widetilde{a}_{ik})= \alpha_{ik}$, $Pc( \widetilde{b}_{jq})= \beta_{jq}$,
$T(\widetilde{a}_{ik})= a_{ik}$, $T(\widetilde{b}_{jq})= b_{jq}$, $f=
``\sum_{i,k} a_{ik} x^i y^k"$, $g= ``\sum_{j,q} b_{jq} x^j y^q"$. Let $I$, $J$
be the support of $f$ and $g$ with respect to $x$. Consider both resultants
\[R(I,J,\mathbb{K})(\widetilde{f},\widetilde{g})= \sum_{r=0}^N \widetilde{h}_r
y^r \textrm{\ and\ } R_t(I,J,\mathbb{K})(f,g)=``\sum_{r=0}^Nh_ry^r".\] It
happens that $T(\widetilde{h}_r)\leq h_r$ and the equality holds if and only if
the term $\gamma_r(\alpha,\beta) t^{{-h_r}}$ of $\widetilde{h}_r$ is different
from $0$. As in the generic case the resultant projects correctly by
Lemma~\ref{tropicalizar_resultante_no_homogenea_char0p}, the polynomials
$\gamma_r$ corresponding to vertices of the subdivision of the Newton polytope
of the resultant polynomial (that in this case is a segment) are nonzero
polynomials in $k[\alpha_{ik},\beta_{jq}]$. If no one of them vanish, the
resultant tropicalizes correctly.
\end{proof}

With all these results we are ready to prove our main result, we can provide a
bijection between the stable intersection of two tropical curves and the
intersection of two generic lifts of the curves. Moreover, sufficient residual
conditions for the genericity can be explicitly computed.

\begin{theorem}\label{interestable}
Let $\widetilde{f}$, $\widetilde{g}\in\mathbb{K}[x,y]$. Then, it can be computed
a finite set of polynomials in the principal coefficients of $\widetilde{f}$,
$\widetilde{g}$ depending only on their tropicalization $f$, $g$ such that, if
no one of them vanish, the tropicalization of the intersection of
$\widetilde{f}$, $\widetilde{g}$ is exactly the stable intersection of $f$ and
$g$. Moreover, the multiplicities are conserved.
\[\sum_{\substack{\widetilde{q}\in \widetilde{f}\cap \widetilde{g}\\
T(\widetilde{q})=q}} \textrm{mult} (\widetilde{q}) =\textrm{mult}_t(q)\]
\end{theorem}
\begin{proof}
Proposition~\ref{resultantes} provides a set $S$ of polynomials in the principal
coefficients of $\widetilde{f}$ and $\widetilde{g}$ such that, if no one
vanishes, the algebraic resultants $\textrm{Res}_x(\widetilde{f},
\widetilde{g})$ and $\textrm{Res}_y(\widetilde{f},\widetilde{g})$ define the
same tropical varieties as $\textrm{Res}_x(f, g)$ and $\textrm{Res}_y(f, g)$.
These two resultants define a finite set $P$ that contains the stable
intersection. The problem is that, in the tropical case, it is possible that the
intersection of $P$ with both curves may be strictly larger than the stable
intersection of the curves, see Example~\ref{resultex}. So, we need another
polynomial in order to discriminate the points in this intersection that are not
stable points. Take $a$, any natural number such that the affine function $x-ay$
is injective in the finite set $P$. Make the monomial change of coordinates
$z=xy^{{-a}}$. The polynomial $\textrm{Res}_y (\widetilde{f}(zy^{a}, y),
\widetilde{g}(zy^{a}, y))= \widetilde{R}(z)= \widetilde{R}(xy^{-a})$ encodes the
values $xy^{{-a}}$ of the common roots of $\widetilde{f}$ and $\widetilde{g}$.
We add to the set $S$ the restrictions in the principal coefficients of this
resultant to be compatible with tropicalization according to
Proposition~\ref{resultantes}. These values $xy^{{-a}}$ of the algebraic
intersection points correspond with the possible values $x-ay$ of the
tropicalization of the roots. As the linear function is injective in $P$, then
$\mathcal{T}(f) \cap \mathcal{T}(g) \cap \mathcal{T}(\textrm{Res}_x(f,g)) \cap
\mathcal{T}(\textrm{Res}_y(f,g)) \cap \mathcal{T}(R(``xy^{{-a}}"))$ is exactly
the tropicalization of the intersection points of any system
$(\widetilde{f},\widetilde{g})$ verifying the restrictions of $S$. By,
Lemma~\ref{generic-lift_stable-int}, this set is contained in the stable
intersection of $f$ and $g$.

To prove that the multiplicities are conserved, consider the field
$\mathbb{K}=\mathbb{C}{((t^\mathbb{R}))}$ of generalized Puiseux series, in this
case \[\sum_{\substack{\widetilde{q} \in \widetilde{f} \cap \widetilde{g}\\
T(\widetilde{q})= q}} \textrm{mult}( \widetilde{q}) \leq \textrm{mult}_t(q).\]
because the sum on the left is bounded by the mixed volume of the residual
polynomials $\widetilde{f}_q$, $\widetilde{g}_q$ over $q$ by
Bernstein-Koushnirenko Theorem (c.f. \cite{bernstein} \cite{koushnirenko}
\cite{rojas}). This mixed volume is, by definition, the tropical multiplicity of
$q$ on the right. On the other hand, the sum on the left is, over any field, the
sum of the multiplicities of the algebraic roots of $\widetilde{R}(xy^{{-a}})$
projecting onto $q$. By the previous results on the correct projection of the
resultant, this multiplicity does not depend on $\mathbb{K}$, because it is the
degree minus the order of the residual polynomial $R(xy^{-a})_{q_x-aq_y}$, or,
equivalently, the multiplicity of $q$ as a root of $T(R(xy^{{-a}}))$. Moreover,
this multiplicity is the mixed volume of the residual polynomials over $q$. That
is, the inequality \[\sum_{\substack{\widetilde{q} \in \widetilde{f} \cap
\widetilde{g}\\ T(\widetilde{q})= q}} \textrm{mult}( \widetilde{q}) \leq
\textrm{mult}_t(q)\] holds for any field. The total number of roots of
$\widetilde{f}$ and $\widetilde{g}$ counted with multiplicities in the torus
equals the sum of multiplicities of the stable roots of $f$ and $g$, because, in
both cases, this is the degree minus the order of $R(xy^{{-a}})$. From this, we
conclude that \[\sum_{\substack{\widetilde{q}\in \widetilde{f}\cap
\widetilde{g}\\
T(\widetilde{q})=q}}\textrm{mult}(\widetilde{q})=\textrm{mult}_t(q)\] Hence, the
projection of the intersection of $\widetilde{f}$ and $\widetilde{g}$ is exactly
the stable intersection.
\end{proof}

Along the proof of the Theorem we have proved the following result, that asserts
that the tropical resultant of two tropical curves has a geometric meaning
analogous to the algebraic resultant.

\begin{corollary}
Let $f$, $g\in\mathbb{T}[x,y]$ be two tropical polynomials. Let $h(y)\in
\mathbb{T}[y]$ be a tropical resultant of $f$ and $g$ with respect to the $x$
variable. Then, the tropical roots of $h$ are exactly the $y$-th coordinates of
the stable intersection of $f$ and $g$.
\end{corollary}

\begin{figure}
\begin{center}
\includegraphics[width=0.5\linewidth]{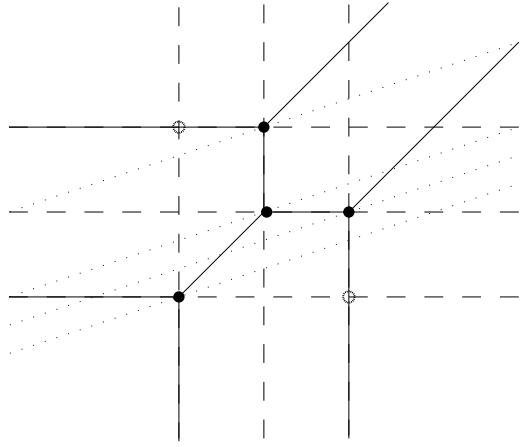}
\caption{Three resultants are needed to compute the stable
intersection.} \label{fig:discrimine}
\end{center}
\end{figure}

\begin{example}\label{resultex}
Consider $f=g=``0+1x+1y+1xy+0x^2+0y^2"$, two conics. Their stable intersection
is the set $\{(-1,-1)$, $(0,1)$, $(1,0)$, $(0,0)\}$. Compute the resultants:
$\textrm{Res}_x(f,g)=``0+1y+1y^2+1y^3+0y^4",$ by symmetry $\textrm{Res}_y(f,
g)=``0+1x+1x^2+1x^3+0x^4"$. Their roots are the lines $y=-1$, $y=0$, $y=1$ and
$x=-1$, $x=0$, $x=1$ respectively. In both cases the multiplicity of the roots
$-1$ and $1$ is 1, while the multiplicity of $0$ is $2$. The intersection of
this lines and the two curves gives the four stable points plus $(-1,1)$ and
$(1,-1)$. We need another resultant that discriminates the points. See
Figure~\ref{fig:discrimine}. Take $x-3y$, the first affine function $x-ay$ that
is injective over these points. $f(``zy^3", y)=``0 +1y +0y^2 +1y^3z +1y^4z
+0y^6z^2"$. $\textrm{Res}_y(f(``zy^3", y), g(``zy^3", y))=``6z^8 +9z^9 +9z^{10}
+8z^{11} +6z^{12}"$. Its roots are $0,1,2,-3$, all with multiplicity 1. It is
easy to check now that the intersection of the two curves and the three
resultants is exactly the stable intersection. The two extra points take the
values -4, 4 in the monomial $``xy^{-3}"$, moreover, every point has
intersection multiplicity equal to one.

Two generic lifts of the cubics are of the form:
\[\widetilde{f} = a_1 +a_xt^{-1}x +a_yt^{-1}y +a_{xy}t^{-1}xy +a_{xx}x^2
+a_{yy}y^2\]
\[\widetilde{g} = c_1 +c_xt^{-1}x +c_yt^{-1}y +c_{xy}t^{-1}xy +c_{xx}x^2
+c_{yy}y^2\]

The residual conditions for the compatibility of the algebraic and tropical
resultant with respect to $x$ are:\\\mbox{}

\noindent{\small
$-\gamma_{xy}$ $\gamma_{xx}$ $\alpha_{xy}$ $\alpha_{yy}$ $-\gamma_{xy}$
$\alpha_{xy}$ $\alpha_{xx}$ $\gamma_{yy}$ $+\gamma_{xy}^2$ $\alpha_{xx}$
$\alpha_{yy}$ $+\gamma_{yy} $ $\gamma_{xx}$ $\alpha_{xy}^2$, $-\gamma_x$
$\gamma_{xx}$ $\alpha_x$ $\alpha_1$ $-\gamma_x$ $\alpha_x$ $\alpha_{xx}$
$\gamma_1$ $+\gamma_1$ $\gamma_{xx}$ $\alpha_{x}^2$ $+\alpha_{xx}$
$\gamma_{x}^2$ $\alpha_{1}$, $\gamma_{y}$ $\gamma_{xx}$ $\alpha_{x}^2$
$-\gamma_{x}$ $\gamma_{xx}$ $\alpha_{x}$ $\alpha_{y}$ $+\alpha_{xx}$
$\gamma_{x}^2$ $\alpha_{y}$ $-\gamma_{x}$ $\alpha_{x}$ $\alpha_{xx}$
$\gamma_{y}$, $-\gamma_{xy}$ $\alpha_{xy}$ $\alpha_{xx}$ $\gamma_{y}$
$+\gamma_{y}$ $\gamma_{xx}$ $\alpha_{xy}^2$ $-\gamma_{xy}$ $\gamma_{xx}$
$\alpha_{xy}$ $\alpha_{y}$ $+\gamma_{xy}^2$ $\alpha_{ xx}$
$\alpha_{y}$}\\\mbox{}

For the resultant with respect to $y$, the compatibility conditions
are:\\\mbox{}

\noindent{\small
$-\gamma_{y}$ $\gamma_{yy}$ $\alpha_{y}$ $\alpha_{1}$ $-\gamma_{y}$ $\alpha_{y}$
$\alpha_{yy}$ $\gamma_{1}$ $+\gamma_{1}$ $\gamma_{yy}$ $\alpha_{y}^2$
$+\gamma_{y}^2$ $\alpha_{yy}$ $\alpha_{1}$, $\gamma_{x}$ $\gamma_{yy}$
$\alpha_{y}^2$ $-\gamma_{y}$ $\alpha_{y}$ $\alpha_{yy}$ $\gamma_{x}$
$+\gamma_{y}^2$ $\alpha_{yy}$ $\alpha_{x}$ $-\gamma_{y}$ $\gamma_{yy}$
$\alpha_{y}$ $\alpha_{x}$, $\gamma_{xy}^2$ $\alpha_{yy}$ $\alpha_{x}$
$+\gamma_{x}$ $\gamma_{yy}$ $\alpha_{xy}^2$ $-\gamma_{xy}$ $\gamma_{yy}$
$\alpha_{xy}$ $\alpha_{x}$ $-\gamma_{xy}$ $\alpha_{xy}$ $\alpha_{yy}$
$\gamma_{x}$, $-\gamma_{xy}$ $\gamma_{xx}$ $\alpha_{xy}$ $\alpha_{yy}$
$-\gamma_{xy}$ $\alpha_{xy}$ $\alpha_{xx}$ $\gamma_{yy}$ $+\gamma_{xy}^2$
$\alpha_{xx}$ $\alpha_{yy}$ $+\gamma_{yy}$ $\gamma_{xx}$
$\alpha_{xy}^2$.}\\\mbox{}

Finally, the third resultant is a degree twelve polynomial in the variable $z$.
The residual conditions for its compatibility with the tropical resultant
are:\\\mbox{}

\noindent{\small
$2\gamma_{yy}^2$ $\gamma_{xx}$ $\alpha_{xy}^3$ $\alpha_{yy}$ $\gamma_{y}$
$\alpha_{y}$ $\gamma_1$ $\alpha_{xx}$ $\gamma_{xy}$ $-\gamma_{yy}^2$
$\gamma_{xx}^2$ $\alpha_{xy}^4$ $\alpha_{yy}$ $\gamma_{y}$ $\alpha_{y}$
$\gamma_1$ $-2\gamma_{yy}^2$ $\alpha_{xy}$ $\gamma_{xy}^3$ $\alpha_{xx}^2$
$\alpha_{y}^2$ $\gamma_1$ $\alpha_{yy}$ $+\gamma_{xy}^4$ $\alpha_{xx}^2$
$\gamma_{yy}$ $\alpha_{yy}^2$ $\alpha_{y}^2$ $\gamma_1$ $-\gamma_{xy}^4$
$\alpha_{xx}^2$ $\gamma_{yy}$ $\alpha_{yy}^2$ $\alpha_{y}$ $\gamma_{y}$
$\alpha_{1}$ $+\gamma_{yy}^2$ $\alpha_{xy}^2$ $\alpha_{yy}$ $\gamma_{y}^2$
$\alpha_{xx}^2$ $\gamma_{xy}^2$ $\alpha_{1}$ $-\gamma_{xy}^2$ $\gamma_{xx}^2$
$\alpha_{xy}^2$ $\alpha_{yy}^3$ $\gamma_{y}$ $\alpha_{y}$ $\gamma_1$
$-2\gamma_{yy}$ $\alpha_{xy}$ $\gamma_{xy}^3$ $\alpha_{xx}^2$ $\alpha_{yy}^2$
$\gamma_{y}^2$ $\alpha_{1}$ $+2\gamma_{yy}^2$ $\gamma_{xx}^2$ $\alpha_{xy}^3$
$\gamma_{xy}$ $\alpha_{y}$ $\gamma_{y}$ $\alpha_{yy}$ $\alpha_{1}$
$-2\gamma_{yy}^2$ $\gamma_{xx}^2$ $\alpha_{xy}^3$ $\gamma_{xy}$ $\alpha_{y}^2$
$\gamma_1$ $\alpha_{yy}$ $-\gamma_{xy}^2$ $\gamma_{xx}^2$ $\alpha_{xy}^2$
$\gamma_{yy}$ $\alpha_{yy}^2$ $\alpha_{y}$ $\gamma_{y}$ $\alpha_{1}$
$+\gamma_{xy}^2$ $\gamma_{xx}^2$ $\alpha_{xy}^2$ $\gamma_{yy}$ $\alpha_{yy}^2$
$\alpha_{y}^2$ $\gamma_1$ $-4\gamma_{yy}^2$ $\gamma_{xx}$ $\alpha_{xy}^2$
$\gamma_{xy}^2$ $\alpha_{xx}$ $\alpha_{y}$ $\gamma_{y}$ $\alpha_{yy}$
$\alpha_{1}$ $-2\gamma_{yy}$ $\gamma_{xx}^2$ $\alpha_{xy}^3$ $\gamma_{xy}$
$\alpha_{yy}^2$ $\gamma_{y}^2$ $\alpha_{1}$ $+2\gamma_{yy}^2$ $\alpha_{xy}$
$\gamma_{xy}^3$ $\alpha_{xx}^2$ $\alpha_{y}$ $\gamma_{y}$ $\alpha_{yy}$
$\alpha_{1}$ $+2\gamma_{yy}$ $\gamma_{xx}^2$ $\alpha_{xy}^3$ $\gamma_{xy}$
$\alpha_{yy}^2$ $\gamma_{y}$ $\alpha_{y}$ $\gamma_1$ $+4\gamma_{yy}$
$\gamma_{xx}$ $\alpha_{xy}^2$ $\gamma_{xy}^2$ $\alpha_{yy}^2$ $\gamma_{y}^2$
$\alpha_{xx}$ $\alpha_{1}$ $+\gamma_{yy}^3$ $\gamma_{xx}^2$ $\alpha_{xy}^4$
$\alpha_{y}^2$ $\gamma_1$ $-4\gamma_{yy}$ $\gamma_{xx}$ $\alpha_{xy}^2$
$\gamma_{xy}^2$ $\alpha_{yy}^2$ $\gamma_{y}$ $\alpha_{xx}$ $\alpha_{y}$
$\gamma_1$ $-\gamma_{yy}^3$ $\gamma_{xx}^2$ $\alpha_{xy}^4$ $\alpha_{y}$
$\gamma_{y}$ $\alpha_{1}$ $+2\gamma_{yy}^3$ $\gamma_{xx}$ $\alpha_{xy}^3$
$\alpha_{y}$ $\gamma_{y}$ $\alpha_{xx}$ $\gamma_{xy}$ $\alpha_{1}$
$-2\gamma_{yy}^3$ $\gamma_{xx}$ $\alpha_{xy}^3$ $\alpha_{y}^2$ $\gamma_1$
$\alpha_{xx}$ $\gamma_{xy}$ $-\gamma_{yy}^3$ $\alpha_{xy}^2$ $\alpha_{xx}^2$
$\gamma_{xy}^2$ $\alpha_{y}$ $\gamma_{y}$ $\alpha_{1}$ $+\gamma_{yy}^3$
$\alpha_{xy}^2$ $\alpha_{xx}^2$ $\gamma_{xy}^2$ $\alpha_{y}^2$ $\gamma_1$
$+\gamma_{xy}^2$ $\gamma_{xx}^2$ $\alpha_{xy}^2$ $\alpha_{yy}^3$ $\gamma_{y}^2$
$\alpha_{1}$ $-\gamma_{yy}^2$ $\alpha_{xy}^2$ $\alpha_{yy}$ $\gamma_{y}$
$\alpha_{xx}^2$ $\gamma_{xy}^2$ $\alpha_{y}$ $\gamma_1$ $-2\gamma_{yy}^2$
$\gamma_{xx}$ $\alpha_{xy}^3$ $\alpha_{xx}$ $\gamma_{xy}$ $\alpha_{yy}$
$\gamma_{y}^2$ $\alpha_{1}$ $+\gamma_{yy}^2$ $\gamma_{xx}^2$ $\alpha_{xy}^4$
$\alpha_{yy}$ $\gamma_{y}^2$ $\alpha_{1}$ $-\gamma_{xy}^4$ $\alpha_{xx}^2$
$\alpha_{yy}^3$ $\gamma_{y}$ $\alpha_{y}$ $\gamma_1$ $+4\gamma_{yy}^2$
$\gamma_{xx}$ $\alpha_{xy}^2$ $\gamma_{xy}^2$ $\alpha_{xx}$ $\alpha_{y}^2$
$\gamma_1$ $\alpha_{yy}$ $+\gamma_{xy}^4$ $\alpha_{xx}^2$ $\alpha_{yy}^3$
$\gamma_{y}^2$ $\alpha_{1}$ $+2\gamma_{xy}^3$ $\alpha_{xx}$ $\gamma_{xx}$
$\alpha_{xy}$ $\gamma_{yy}$ $\alpha_{yy}^2$ $\alpha_{y}$ $\gamma_{y}$
$\alpha_{1}$ $-2\gamma_{xy}^3$ $\alpha_{xx}$ $\gamma_{xx}$ $\alpha_{xy}$
$\gamma_{yy}$ $\alpha_{yy}^2$ $\alpha_{y}^2$ $\gamma_1$ $-2\gamma_{xy}^3$
$\gamma_{xx}$ $\alpha_{xy}$ $\alpha_{xx}$ $\alpha_{yy}^3$ $\gamma_{y}^2$
$\alpha_{1}$ $+2\gamma_{xy}^3$ $\gamma_{xx}$ $\alpha_{xy}$ $\alpha_{xx}$
$\alpha_{yy}^3$ $\gamma_{y}$ $\alpha_{y}$ $\gamma_1$ $+2\gamma_{yy}$
$\alpha_{xy}$ $\gamma_{xy}^3$ $\alpha_{xx}^2$ $\alpha_{yy}^2$ $\gamma_{y}$
$\alpha_{y}$ $\gamma_1,$\\
$3\gamma_{xy}$ $\gamma_{xx}^2$ $\alpha_{xy}^4$ $\gamma_{y}^2$ $\alpha_{y}^2$
$\gamma_1$ $-3\gamma_{xy}$ $\gamma_{xx}^2$ $\alpha_{xy}^4$ $\gamma_{y}^3$
$\alpha_{y}$ $\alpha_{1}$ $-\gamma_{xx}^2$ $\alpha_{xy}^5$ $\gamma_{y}^3$
$\alpha_{y}$ $\gamma_1$ $+3\gamma_{xy}^3$ $\alpha_{xx}^2$ $\gamma_{y}^2$
$\alpha_{xy}^2$ $\alpha_{y}^2$ $\gamma_1$ $-\gamma_{xy}^5$ $\alpha_{xx}^2$
$\alpha_{y}^3$ $\gamma_{y}$ $\alpha_{1}$ $+\gamma_{xy}^3$ $\gamma_{xx}^2$
$\alpha_{xy}^2$ $\alpha_{y}^4$ $\gamma_1$ $+6\gamma_{xy}^3$ $\alpha_{xx}$
$\gamma_{y}$ $\gamma_{xx}$ $\alpha_{xy}^2$ $\alpha_{y}^3$ $\gamma_1$
$-3\gamma_{xy}^4$ $\alpha_{xx}^2$ $\gamma_{y}$ $\alpha_{xy}$ $\alpha_{y}^3$
$\gamma_1$ $-6\gamma_{xy}^3$ $\alpha_{xx}$ $\gamma_{y}^2$ $\gamma_{xx}$
$\alpha_{xy}^2$ $\alpha_{y}^2$ $\alpha_{1}$ $+3\gamma_{xy}^4$ $\alpha_{xx}^2$
$\gamma_{y}^2$ $\alpha_{xy}$ $\alpha_{y}^2$ $\alpha_{1}$ $+\gamma_{xy}^5$
$\alpha_{xx}^2$ $\alpha_{y}^4$ $\gamma_1$ $-3\gamma_{xy}^3$ $\alpha_{xx}^2$
$\gamma_{y}^3$ $\alpha_{xy}^2$ $\alpha_{y}$ $\alpha_{1}$ $-2\gamma_{xy}^4$
$\gamma_{xx}$ $\alpha_{xy}$ $\alpha_{xx}$ $\alpha_{y}^4$ $\gamma_1$
$+2\gamma_{xy}$ $\gamma_{xx}$ $\alpha_{xy}^4$ $\gamma_{y}^3$ $\alpha_{xx}$
$\alpha_{y}$ $\gamma_1$ $-\gamma_{xy}^2$ $\alpha_{xx}^2$ $\gamma_{y}^3$
$\alpha_{xy}^3$ $\alpha_{y}$ $\gamma_1$ $-2\gamma_{xy}$ $\gamma_{xx}$
$\alpha_{xy}^4$ $\gamma_{y}^4$ $\alpha_{xx}$ $\alpha_{1}$ $+2\gamma_{xy}^4$
$\gamma_{xx}$ $\alpha_{xy}$ $\alpha_{xx}$ $\alpha_{y}^3$ $\gamma_{y}$
$\alpha_{1}$ $-\gamma_{xy}^3$ $\gamma_{xx}^2$ $\alpha_{xy}^2$ $\alpha_{y}^3$
$\gamma_{y}$ $\alpha_{1}$ $+\gamma_{xx}^2$ $\alpha_{xy}^5$ $\gamma_{y}^4$
$\alpha_{1}$ $+3\gamma_{xy}^2$ $\gamma_{xx}^2$ $\alpha_{xy}^3$ $\gamma_{y}^2$
$\alpha_{y}^2$ $\alpha_{1}$ $-6\gamma_{xy}^2$ $\alpha_{xx}$ $\gamma_{y}^2$
$\gamma_{xx}$ $\alpha_{xy}^3$ $\alpha_{y}^2$ $\gamma_1$ $+6\gamma_{xy}^2$
$\alpha_{xx}$ $\gamma_{y}^3$ $\gamma_{xx}$ $\alpha_{xy}^3$ $\alpha_{y}$
$\alpha_{1}$ $-3\gamma_{xy}^2$ $\gamma_{xx}^2$ $\alpha_{xy}^3$ $\gamma_{y}$
$\alpha_{y}^3$ $\gamma_1$ $+\gamma_{xy}^2$ $\alpha_{xx}^2$ $\gamma_{y}^4$
$\alpha_{xy}^3$ $\alpha_{1},$\\
$\gamma_{xy}^3$ $\alpha_{xx}^2$ $\gamma_{x}^2$ $\alpha_{x}$ $\alpha_{y}^3$
$\gamma_1$ $+\gamma_{x}$ $\gamma_{xx}^2$ $\alpha_{xy}^3$ $\alpha_{x}^2$
$\gamma_{y}^3$ $\alpha_{1}$ $+\gamma_{x}^3$ $\gamma_{xx}^2$ $\alpha_{xy}^3$
$\alpha_{y}^2$ $\gamma_{y}$ $\alpha_{1}$ $+\gamma_{xy}^3$ $\alpha_{xx}^2$
$\alpha_{x}^3$ $\gamma_{y}^2$ $\alpha_{y}$ $\gamma_1$ $-\gamma_{x}^3$
$\alpha_{xy}$ $\alpha_{xx}^2$ $\gamma_{xy}^2$ $\alpha_{y}^3$ $\gamma_1$
$+2\gamma_{xy}$ $\gamma_{xx}^2$ $\alpha_{xy}^2$ $\gamma_{x}$ $\alpha_{x}^2$
$\gamma_{y}^2$ $\alpha_{y}$ $\alpha_{1}$ $+2\gamma_{xy}^2$ $\alpha_{xx}$
$\gamma_{xx}$ $\alpha_{x}^3$ $\gamma_{y}^3$ $\alpha_{xy}$ $\alpha_{1}$
$+4\gamma_{xy}^2$ $\alpha_{xx}$ $\gamma_{x}$ $\gamma_{xx}$ $\alpha_{x}^2$
$\gamma_{y}$ $\alpha_{xy}$ $\alpha_{y}^2$ $\gamma_1$ $-4\gamma_{xy}^2$
$\alpha_{xx}$ $\gamma_{x}$ $\gamma_{xx}$ $\alpha_{x}^2$ $\gamma_{y}^2$
$\alpha_{xy}$ $\alpha_{y}$ $\alpha_{1}$ $-2\gamma_{xy}^3$ $\alpha_{xx}^2$
$\gamma_{x}$ $\alpha_{x}^2$ $\alpha_{y}^2$ $\gamma_{y}$ $\gamma_1$
$+2\gamma_{xy}^3$ $\alpha_{xx}^2$ $\gamma_{x}$ $\alpha_{x}^2$ $\alpha_{y}$
$\gamma_{y}^2$ $\alpha_{1}$ $-\gamma_{xy}^3$ $\alpha_{xx}^2$ $\gamma_{x}^2$
$\alpha_{x}$ $\alpha_{y}^2$ $\gamma_{y}$ $\alpha_{1}$ $-\gamma_{x}$
$\gamma_{xx}^2$ $\alpha_{xy}^3$ $\alpha_{x}^2$ $\gamma_{y}^2$ $\alpha_{y}$
$\gamma_1$ $+\gamma_{xy}$ $\gamma_{xx}^2$ $\alpha_{xy}^2$ $\gamma_{x}^2$
$\alpha_{x}$ $\alpha_{y}^3$ $\gamma_1$ $-\gamma_{xy}$ $\gamma_{xx}^2$
$\alpha_{xy}^2$ $\gamma_{x}^2$ $\alpha_{x}$ $\alpha_{y}^2$ $\gamma_{y}$
$\alpha_{1}$ $+2\gamma_{x}^2$ $\gamma_{xx}^2$ $\alpha_{xy}^3$ $\gamma_{y}$
$\alpha_{y}^2$ $\gamma_1$ $\alpha_{x}$ $-2\gamma_{x}^2$ $\gamma_{xx}^2$
$\alpha_{xy}^3$ $\gamma_{y}^2$ $\alpha_{y}$ $\alpha_{x}$ $\alpha_{1}$
$+2\gamma_{x}^3$ $\gamma_{xx}$ $\alpha_{xy}^2$ $\gamma_{xy}$ $\alpha_{y}^3$
$\gamma_1$ $\alpha_{xx}$ $-2\gamma_{x}^3$ $\gamma_{xx}$ $\alpha_{xy}^2$
$\gamma_{xy}$ $\alpha_{y}^2$ $\gamma_{y}$ $\alpha_{xx}$ $\alpha_{1}$
$+\gamma_{xy}$ $\gamma_{xx}^2$ $\alpha_{xy}^2$ $\alpha_{x}^3$ $\gamma_{y}^2$
$\alpha_{y}$ $\gamma_1$ $+\gamma_{x}^3$ $\alpha_{xy}$ $\alpha_{xx}^2$
$\gamma_{xy}^2$ $\alpha_{y}^2$ $\gamma_{y}$ $\alpha_{1}$ $+2\gamma_{x}$
$\gamma_{xx}$ $\alpha_{xy}^2$ $\gamma_{xy}$ $\alpha_{x}^2$ $\gamma_{y}^2$
$\alpha_{xx}$ $\alpha_{y}$ $\gamma_1$ $-2\gamma_{x}$ $\gamma_{xx}$
$\alpha_{xy}^2$ $\gamma_{xy}$ $\alpha_{x}^2$ $\gamma_{y}^3$ $\alpha_{xx}$
$\alpha_{1}$ $-4\gamma_{x}^2$ $\gamma_{xx}$ $\alpha_{xy}^2$ $\gamma_{xy}$
$\alpha_{x}$ $\gamma_{y}$ $\alpha_{y}^2$ $\gamma_1$ $\alpha_{xx}$
$+4\gamma_{x}^2$ $\gamma_{xx}$ $\alpha_{xy}^2$ $\gamma_{xy}$ $\alpha_{x}$
$\gamma_{y}^2$ $\alpha_{y}$ $\alpha_{xx}$ $\alpha_{1}$ $+2\gamma_{x}^2$
$\alpha_{xy}$ $\gamma_{xy}^2$ $\alpha_{x}$ $\alpha_{y}^2$ $\alpha_{xx}^2$
$\gamma_{y}$ $\gamma_1$ $-2\gamma_{x}^2$ $\alpha_{xy}$ $\gamma_{xy}^2$
$\alpha_{x}$ $\alpha_{y}$ $\alpha_{xx}^2$ $\gamma_{y}^2$ $\alpha_{1}$
$-2\gamma_{xy}^2$ $\alpha_{xx}$ $\gamma_{x}^2$ $\gamma_{xx}$ $\alpha_{x}$
$\alpha_{xy}$ $\alpha_{y}^3$ $\gamma_1$ $+2\gamma_{xy}^2$ $\alpha_{xx}$
$\gamma_{x}^2$ $\gamma_{xx}$ $\alpha_{x}$ $\alpha_{xy}$ $\alpha_{y}^2$
$\gamma_{y}$ $\alpha_{1}$ $-\gamma_{x}^3$ $\gamma_{xx}^2$ $\alpha_{xy}^3$
$\alpha_{y}^3$ $\gamma_1$ $-\gamma_{xy}^3$ $\alpha_{xx}^2$ $\alpha_{x}^3$
$\gamma_{y}^3$ $\alpha_{1}$ $-2\gamma_{xy}^2$ $\alpha_{xx}$ $\gamma_{xx}$
$\alpha_{x}^3$ $\gamma_{y}^2$ $\alpha_{y}$ $\gamma_1$ $\alpha_{xy}$
$-\gamma_{xy}$ $\gamma_{xx}^2$ $\alpha_{xy}^2$ $\alpha_{x}^3$ $\gamma_{y}^3$
$\alpha_{1}$ $-2\gamma_{xy}$ $\gamma_{xx}^2$ $\alpha_{xy}^2$ $\gamma_{x}$
$\alpha_{x}^2$ $\gamma_{y}$ $\alpha_{y}^2$ $\gamma_1$ $-\gamma_{x}$
$\alpha_{xx}^2$ $\gamma_{xy}^2$ $\alpha_{x}^2$ $\gamma_{y}^2$ $\alpha_{xy}$
$\alpha_{y}$ $\gamma_1$ $+\gamma_{x}$ $\alpha_{xx}^2$ $\gamma_{xy}^2$
$\alpha_{x}^2$ $\gamma_{y}^3$ $\alpha_{xy}$ $\alpha_{1},$\\
$6$ $\gamma_{xx}^2$ $\alpha_{xx}$ $\gamma_{x}^2$ $\alpha_{x}^3$ $\gamma_{y}$
$\alpha_{y}^2$ $\gamma_1$ $-\gamma_{xx}^3$ $\alpha_{x}^3$ $\gamma_{x}^2$
$\alpha_{y}^2$ $\gamma_{y}$ $\alpha_{1}$ $-\gamma_{xx}^3$ $\alpha_{x}^5$
$\gamma_{y}^3$ $\alpha_{1}$ $-6\gamma_{xx}^2$ $\alpha_{xx}$ $\gamma_{x}^2$
$\alpha_{x}^3$ $\gamma_{y}^2$ $\alpha_{y}$ $\alpha_{1}$ $+6\gamma_{xx}$
$\alpha_{xx}^2$ $\gamma_{x}^3$ $\gamma_{y}^2$ $\alpha_{x}^2$ $\alpha_{y}$
$\alpha_{1}$ $-\alpha_{xx}^3$ $\gamma_{x}^5$ $\alpha_{y}^3$ $\gamma_1$
$+\gamma_{xx}^3$ $\alpha_{x}^5$ $\gamma_{y}^2$ $\alpha_{y}$ $\gamma_1$
$+3\gamma_{xx}^2$ $\alpha_{xx}$ $\gamma_{x}^3$ $\alpha_{x}^2$ $\alpha_{y}^2$
$\gamma_{y}$ $\alpha_{1}$ $+\gamma_{xx}^3$ $\alpha_{x}^3$ $\gamma_{x}^2$
$\alpha_{y}^3$ $\gamma_1$ $+\alpha_{xx}^3$ $\gamma_{x}^5$ $\alpha_{y}^2$
$\gamma_{y}$ $\alpha_{1}$ $-\alpha_{xx}^3$ $\gamma_{x}^3$ $\alpha_{x}^2$
$\gamma_{y}^2$ $\alpha_{y}$ $\gamma_1$ $+3\gamma_{xx}^2$ $\alpha_{xx}$
$\gamma_{x}$ $\alpha_{x}^4$ $\gamma_{y}^3$ $\alpha_{1}$ $-6\gamma_{xx}$
$\alpha_{xx}^2$ $\gamma_{x}^3$ $\gamma_{y}$ $\alpha_{x}^2$ $\alpha_{y}^2$
$\gamma_1$ $+2\gamma_{xx}^3$ $\alpha_{x}^4$ $\gamma_{x}$ $\gamma_{y}^2$
$\alpha_{y}$ $\alpha_{1}$ $-2\gamma_{xx}^3$ $\alpha_{x}^4$ $\gamma_{x}$
$\gamma_{y}$ $\alpha_{y}^2$ $\gamma_1$ $+\alpha_{xx}^3$ $\gamma_{x}^3$
$\alpha_{x}^2$ $\gamma_{y}^3$ $\alpha_{1}$ $-3\gamma_{xx}^2$ $\alpha_{xx}$
$\gamma_{x}^3$ $\alpha_{x}^2$ $\alpha_{y}^3$ $\gamma_1$ $+3\gamma_{xx}$
$\alpha_{xx}^2$ $\gamma_{x}^4$ $\alpha_{x}$ $\alpha_{y}^3$ $\gamma_1$
$-3\gamma_{xx}$ $\alpha_{xx}^2$ $\gamma_{x}^4$ $\alpha_{x}$ $\alpha_{y}^2$
$\gamma_{y}$ $\alpha_{1}$ $-3\gamma_{xx}^2$ $\alpha_{xx}$,$\gamma_{x}$
$\alpha_{x}^4$ $\gamma_{y}^2$ $\alpha_{y}$ $\gamma_1$ $-3\gamma_{xx}$
$\alpha_{xx}^2$ $\gamma_{x}^2$ $\alpha_{x}^3$ $\gamma_{y}^3$ $\alpha_{1}$
$-2\alpha_{xx}^3$ $\gamma_{x}^4$ $\alpha_{x}$ $\gamma_{y}^2$ $\alpha_{y}$
$\alpha_{1}$ $+2\alpha_{xx}^3$ $\gamma_{x}^4$ $\alpha_{x}$ $\gamma_{y}$
$\alpha_{y}^2$ $\gamma_1$ $+3\gamma_{xx}$ $\alpha_{xx}^2$ $\gamma_{x}^2$
$\alpha_{x}^3$ $\gamma_{y}^2$ $\alpha_{y}$ $\gamma_1,$\\
$3\alpha_{xx}^3$ $\gamma_{x}^4$ $\alpha_{x}^2$ $\alpha_{1}$ $\gamma_1^2$
$+3\gamma_{xx}^3$ $\alpha_{x}^4$ $\gamma_{x}^2$ $\gamma_1$ $\alpha_{1}^2$
$+\gamma_{xx}^3$ $\alpha_{x}^6$ $\gamma_1^3$ $+\alpha_{xx}^3$ $\gamma_{x}^6$
$\alpha_{1}^3$ $-3\gamma_{xx}^3$ $\alpha_{x}^5$ $\gamma_{x}$ $\gamma_1^2$
$\alpha_{1}$ $+9\gamma_{xx}$ $\alpha_{xx}^2$ $\gamma_{x}^4$ $\alpha_{x}^2$
$\alpha_{1}^2$ $\gamma_1$ $+3\gamma_{xx}$ $\alpha_{xx}^2$ $\gamma_{x}^2$
$\alpha_{x}^4$ $\gamma_1^3$ $+3\gamma_{xx}^2$ $\alpha_{xx}$ $\gamma_{x}^4$
$\alpha_{x}^2$ $\alpha_{1}^3$ $-3\alpha_{xx}^3$ $\gamma_{x}^5$ $\alpha_{x}$
$\alpha_{1}^2$ $\gamma_1$ $-9\gamma_{xx}^2$ $\alpha_{xx}$ $\gamma_{x}^3$
$\alpha_{x}^3$ $\gamma_1$ $\alpha_{1}^2$ $-3\gamma_{xx}^2$ $\alpha_{xx}$
$\gamma_{x}$ $\alpha_{x}^5$ $\gamma_1^3$ $-3\gamma_{xx}$ $\alpha_{xx}^2$
$\gamma_{x}^5$ $\alpha_{x}$ $\alpha_{1}^3$ $+9\gamma_{xx}^2$ $\alpha_{xx}$
$\gamma_{x}^2$ $\alpha_{x}^4$ $\gamma_1^2$ $\alpha_{1}$ $-\alpha_{xx}^3$
$\gamma_{x}^3$ $\alpha_{x}^3$ $\gamma_1^3$ $-\gamma_{xx}^3$ $\alpha_{x}^3$
$\gamma_{x}^3$ $\alpha_{1}^3$ $-9\gamma_{xx}$ $\alpha_{xx}^2$ $\gamma_{x}^3$
$\alpha_{x}^3$ $\gamma_1^2$ $\alpha_{1}$}
\end{example}

\section{Some Remarks}\label{sec_remarks}

As a consequence of Theorem~\ref{interestable}, a new proof of
Bernstein-Koushnirenko Theorem for plane curves over an arbitrary algebraically
closed field can be derived from the classic Theorem over $\mathbb{C}$
(\cite{bernstein}, \cite{koushnirenko}). We refer to \cite{rojas} for a direct
proof in positive characteristic.

\begin{corollary}
Let $\widetilde{f}$, $\widetilde{g}$ be two polynomials over $\mathbb{K}$, an
algebraically closed. Let $\Delta_f$, $\Delta_g$ be the Newton polytope of the
polynomials $\widetilde{f}$ and $\widetilde{g}$ respectively. Then, if the
coefficients of $\widetilde{f}$ and $\widetilde{g}$ are generic, then the number
of common roots of the curves in $(\mathbb{K}^*)^2$ counted with multiplicities
is the mixed volume of the Newton polygons \[\mathcal{M}(\Delta_f,
\Delta_g)=vol(\Delta_f + \Delta_g) -vol(\Delta_f) - vol(\Delta_g)\]
\end{corollary}
\begin{proof}
If the coefficients of the polynomials are generic, the number of roots in the
torus counted with multiplicities is the degree minus the order of the resultant
of the two polynomials with respect to one of the variables. This number does
only depend on the support of the polynomials, and it is equal to the mixed
volume of the Newton polygons, because this is the number of stable intersection
points of two tropical curves of Newton polygons $\Delta_f$, $\Delta_g$.
\end{proof}

\begin{remark}
Another application of the techniques developed in this article is the
computation of tropical bases. Theorem~\ref{initialserie} proves that for a
hypersurface $\widetilde{f}$, the projection $T(\{\widetilde{f}=0\})=
\mathcal{T}(f)$. This is not true for general ideals. If
$\mathcal{I}=(\widetilde{f}_1, \ldots, \widetilde{f}_m) \subseteq
\mathbb{K}[x_1, \ldots, x_n]$ and $\mathcal{V}$ is the variety it defines in
$(\mathbb{K}^*)^n$, \[T(\mathcal{V}) \subseteq \bigcap_{i=1}^m
\mathcal{T}(f_i),\] but it is possible that both sets are different. A set of
generators $\widetilde{g}_1, \ldots, \widetilde{g}_r$ of $\mathcal{I}$ such that
$T(\mathcal{V}) = \bigcap_{i=1}^r \mathcal{T}(g_r)$ is called a tropical basis
of $\mathcal{I}$. In \cite{Computing_trop_var}, it is proved that every ideal
has a tropical basis and it is provided an algorithm for the case of a prime
ideal $\mathcal{I}$.

An alternative for the computation of a tropical basis of a zero dimensional
ideal in two variables is the following. Let $\mathcal{I}=(\widetilde{f},
\widetilde{g})$ be a zero dimensional ideal in two variables. Let
$\widetilde{R}_x$, $\widetilde{R}_y$ be the resultants with respect to $x$ and
$y$ of the curves. Let $P$ be the intersection of the projections $R_x$ and
$R_y$. This is always a finite set that contains the projection of the
intersection of $\widetilde{f}$, $\widetilde{g}$. It may happen that $P$ is not
contained in the stable intersection of the corresponding tropical curves $f$
and $g$, though. Let $a$ be a natural number such that $x-ay$ is injective in
$P$. Let $\widetilde{R}_z=\textrm{Res}_y (\widetilde{f}(zy^{a},y),
\widetilde{g}(zy^{a},y))$ be another resultant. Then, it follows that
$(\widetilde{f}, \widetilde{g}, \widetilde{R}_x, \widetilde{R}_y,
\widetilde{R}_z)$ is a tropical basis of the ideal $(\widetilde{f},
\widetilde{g})$. This alternative approach is very similar to the regular
projection method that has been developed by Hept and Theobald
\cite{TR_bas_reg_proy}.
\end{remark}

\begin{remark}
Along the article, the notion of tropical resultant has been defined as the
projection of the algebraic resultant. It is needed a precomputation of the
algebraic resultant in order to tropicalize it. For the case of plane curves, it
would be preferable to have a determinantal formula. That is, to prove that the
determinant of the Sylvester matrix of two polynomials define the resultant
variety. But the proof of the properties is achieved by a careful look to the
polynomials involved, paying special attention to the cancellation of terms. In
the case of the determinant of the Sylvester matrix, the tropical determinant of
the Sylvester matrix is the projection of the permanent of the algebraic
determinant. There are cancellation of terms even in the equicharacteristic zero
case. It is conjectured that still the determinant of the Sylvester matrix is a
tropical polynomial that defines the same tropical variety as the resultant
does. The author has checked that it is the case for polynomials up to degree
four with full support.
\end{remark}

\newcommand{\etalchar}[1]{$^{#1}$}

\noindent Luis Felipe Tabera\\
IMDEA Matem\'aticas\\
Facultad de Ciencias C-IX\\
Campus Universidad Aut\'onoma de Madrid, E-28049 Madrid, Spain\\
e-mail : luis.tabera@imdea.org

\begin{thebibliography}{RGST05}
\bibitem[Ber75]{bernstein}
D.~N. Bernstein.
\newblock The number of roots of a system of equations.
\newblock {\em Akademija Nauk SSSR. Funkcional{\/$'$} nyi Analiz i ego Prilo\v
  zenija}, 9(3):1--4, 1975.

\bibitem[BJS{\etalchar{+}}07]{Computing_trop_var}
T.~Bogart, A.~N. Jensen, D.~Speyer, B.~Sturmfels, and R.~R. Thomas.
\newblock Computing tropical varieties.
\newblock {\em J. Symbolic Comput.}, 42(1-2):54--73, 2007.

\bibitem[EKL06]{einsiedler-2004-}
M. Einsiedler, M. Kapranov, and D. Lind.
\newblock Non-{A}rchimedean amoebas and tropical varieties.
\newblock {\em J. Reine Angew. Math.}, 601:139--157, 2006.

\bibitem[GKZ90]{GKZ-polytope-resultant}
I.~M. Gel{$'$}fand, M.~M. Kapranov, and A.~V. Zelevinsky.
\newblock Newton polytopes of the classical resultant and discriminant.
\newblock {\em Advances in Mathematics}, 84(2):237--254, 1990.

\bibitem[HT07]{TR_bas_reg_proy}
K. Hept and T. Theobald.
\newblock Tropical bases by regular projections.
\newblock {\em Preprint}, 2007.
\newblock http://arxiv.org/abs/0708.1727

\bibitem[JMM07]{Lifting-Constr}
A.~N. Jensen, H. Markwig, and T. Markwig.
\newblock An algorithm for lifting points in a tropical variety.
\newblock {\em Preprint}, 2007.
\newblock http://arxiv.org/abs/0705.2441

\bibitem[Kus76]{koushnirenko}
A.~G. Kushnirenko.
\newblock Newton polytopes and the bezout theorem.
\newblock {\em Functional Analysis and Its Applications}, 10(3):233--235, 1976.

\bibitem[Mik05]{Mik05}
G. Mikhalkin.
\newblock Enumerative tropical algebraic geometry in {$\mathbb{R}\sp 2$}.
\newblock {\em Journal of the American Mathematical Society}, 18(2):313--377
  (electronic), 2005.

\bibitem[RGST05]{rgst}
J. Richter-Gebert, B. Sturmfels, and T. Theobald. 
\newblock First steps in tropical geometry.
\newblock In {\em Idempotent mathematics and mathematical physics}, volume 377
  of {\em Contemp. Math.}, pages 289--317. Amer. Math. Soc., Providence, RI,
  2005.

\bibitem[Roj99]{rojas}
J.~M. Rojas.
\newblock Toric intersection theory for affine root counting.
\newblock {\em Journal of Pure and Applied Algebra}, 136(1):67--100, 1999.

\bibitem[Stu94]{Sturmfels-polytope_resultant}
B. Sturmfels.
\newblock On the {N}ewton polytope of the resultant.
\newblock {\em Journal of Algebraic Combinatorics. An International Journal},
  3(2):207--236, 1994.

\bibitem[Stu02]{Stu02}
B. Sturmfels.
\newblock {\em Solving systems of polynomial equations}, volume~97 of {\em CBMS
  Regional Conference Series in Mathematics}.
\newblock Published for the Conference Board of the Mathematical Sciences,
  Washington, DC, 2002.

\bibitem[Tab05]{Pappus-trop}
L.~F. Tabera.
\newblock Tropical constructive {P}appus' theorem.
\newblock {\em International Mathematics Research Notices},
  2005(39):2373--2389, 2005.

\bibitem[Tab06]{Kapranov-EACA}
L.~F. Tabera.
\newblock Constructive proof of extended {K}apranov theorem.
\newblock In {\em Actas del X Encuentro de {\'A}lgebra Computacional y
  Aplicaciones, EACA 2006}, pages 178--181, 2006.
\end{thebibliography}
\end{document}